\documentclass[a4paper,12pt]{article}
\usepackage{times, url}
\usepackage{graphicx,amsmath,amssymb,amsthm,amsfonts,bm,float,cite}
\usepackage{xcolor}
\usepackage[none]{hyphenat}
\usepackage[center]{caption}
\usepackage[pdfencoding=auto]{hyperref}
\usepackage{bookmark}
\usepackage{mathtools}

\begin{document}

\begin{center}
{\LARGE \bf Congruence speed of tetration bases ending with $0$}
\vspace{12mm}

{\Large Marco Rip\`a}

\end{center}
\vspace{16mm}

\noindent {\bf Abstract:} For every non-negative integer $a$ and positive integer $b$, the congruence speed of the tetration $^{b}a$ is the difference between the number of the rightmost digits of $^{b}a$ that are the same as those of $^{b+1}a$ and the number of the rightmost digits of $^{b-1}a$ that are the same as those of $^{b}a$. In the decimal numeral system, if the given base $a$ is not a multiple of $10$, as $b:=b(a)$ becomes sufficiently large, we know that the value of the congruence speed does not depend on $b$ anymore, otherwise the number of the new rightmost zeros of $^{b}a$ drastically increases for any unit increment of $b$ and, for this reason, we have not previously described the congruence speed of $a$ when it is a multiple of $10$. This short note fills the gap by giving the formula for the congruence speed of the mentioned values of $a$ at any given height of the hyperexponent.\\
{\bf Keywords:} Tetration, Exponentiation, Congruence speed, Modular arithmetic, Recurring digits.\\ 
{\bf 2020 Mathematics Subject Classification:} Primary 11A07; Secondary 11D79. 
\vspace{5mm}

%%%%%%%%%%%%%%%%%%%%%%%%%
%%%%%%% SECTION 1 %%%%%%%
%%%%%%%%%%%%%%%%%%%%%%%%%

\section{Introduction} \label{sec:Intr}

For any given pair of non-negative integers $(a, b)$, we call \textit{tetration} the $4$th degree hyperoperator which is recursively defined as
\[^{b}a = \begin{cases} 1 \mbox{ }\mbox{ }\mbox{ }\mbox{ }\mbox{ }\mbox{ }\mbox{ }\mbox{ }\mbox{ \hspace{1.2mm} if $b = 0$} \\
a^{(^{(b-1)}a)}  \mbox{ if $b\geq 1$}
\end{cases}\hspace{-3.5mm}.\]

For this purpose, it is interesting to note that we can also extend the above to a negative value of the hyperexponent since, assuming $\overline{b} \in \mathbb{N}_0 \cup \{-1\}$, $^{\overline{b}}a = \log_{a}\left(^{\overline{b}+1}a \right)$ allows us to state $^{-1}a=\log_{a}\left(a^0 \right)=\log_{a}(1)=0$ and even go beyond (see References \cite{0} and \cite{4}, Section 3).

Assuming radix-$10$ (i.e., the decimal numeral system), for every $a$ that is not a multiple of $10$, it is well known \cite{2} that tetration has a peculiar property, named by the author as ``constancy of the congruence speed'', which we can translate as the fact that, for every sufficiently large $b \coloneqq b(a)$ (for a tight bound, consider the sufficient condition $\tilde{v}(a)+2$, where
$$
\tilde{v}(a):= \begin{cases}v_{5}(a-1) & \textnormal { iff } a \equiv 1\hspace{-3mm}\pmod{5} \\ v_{5}\left(a^{2}+1\right) & \textnormal { iff } a \equiv\{2,3\}\hspace{-3mm}\pmod{5} \\ v_{5}(a+1) & \textnormal { iff } a \equiv 4\hspace{-3mm}\pmod{5} \\ v_{2}\left(a^{2}-1\right)-1 & \textnormal { iff } a \equiv 5\hspace{-3mm}\pmod{10} \end{cases},
$$ while $v_{2}(m)$ and $v_{5}(m)$ indicates the $2$-adic and the $5$-adic valuation of $m$, respectively), the number of the rightmost digits of $^{b}a$ that are the same as those of $^{b+1}a$ minus the number of the rightmost digits of $^{b-1}a$ that remain the same by moving to $^{b}a$ is equal to the number of the new rightmost digits that freeze when we move from $^{b+k}a$ to $^{b+k+1}a$, for every $k \in \mathbb{Z}^+$ (see Definitions~1.1~to~1.3 of Reference \cite{1}).

Now, for each non-negative integer $a$, let us indicate the number of new frozen digits at height $b \in \mathbb{Z}^+$ as $V(a,b)$, while $V(a)$ denotes the (fixed) value of the constant congruence speed of $a$ (i.e., the non-negative integer returned by $V(a,b)$ for any $b \geq \tilde{v}(a)+2$).

As an example, let us take $a=807$. Then, we have $V(807, 1)=0$, $V(807, 2)=4$, $V(807, 3)=4$, $V(807, 4)=4$, $V(807, 5)=4$, and finally $V(807, 6)=V(807, 7)=\cdots=V(807)=3$ since $\tilde{v}_{5}\left(807^2+1\right) + 2= 4+2 = 6$ (by Definition 2.1 of Reference \cite{1}) and

\noindent$^{1}807=807$;

\noindent$^{2}807 \equiv 54962039628331827388850173\color{cyan}{\bf{7943}}$ $\pmod {10^{30}}$;

\noindent$^{3}807 \equiv 6016926514668229405256\color{cyan}{\bf{3285}}\color{red}7943$ $\pmod {10^{30}}$;

\noindent$^{4}807 \equiv 146336906474874632\color{cyan}{\bf{6260}}\color{red}32857943$ $\pmod {10^{30}}$;

\noindent$^{5}807 \equiv 35503490744897\color{cyan}{\bf{3150}}\color{red}626032857943$ $\pmod {10^{30}}$;

\noindent$^{6}807 \equiv 47863568981\color{cyan}{\bf{228}}\color{red}3150626032857943$ $\pmod {10^{30}}$;

\noindent$^{7}807 \equiv 02704888\color{cyan}{\bf{876}}\color{red}2283150626032857943$ $\pmod {10^{30}}$;

\hspace{0.1mm}\vdots

\noindent$^{(1 \textnormal{googolplex})}807 \equiv \color{red} 803001638762283150626032857943$ $\pmod {10^{30}}$.

References \cite{1, 2, 3} describe the value of $V(a,b)$ for any $a$ that is not multiple of $10$.\linebreak In particular, Equation (16) of Reference \cite{1} gives the exact value of $V(a)$ for all the mentioned tetration bases.

Accordingly, in Section \ref{sec:2}, we provide the formula for every $V(a,b)$ such that $a$ is congruent to $0$ modulo $10$ (we point out that if $a$ is a multiple of $10$, then there does not exist any $b$ which stabilizes the congruence speed of the base and thus $V(a)$ is always undefined).

%%%%%%%%%%%%%%%%%%%%%%%%%
%%%%%%% SECTION 2 %%%%%%%
%%%%%%%%%%%%%%%%%%%%%%%%%

\section{Congruence speed of all the tetration bases congruent to \texorpdfstring{$0$ } -modulo \texorpdfstring{$10$}{Congruence speed of all the tetration bases congruent to zero modulo ten}} \label{sec:2}

Let $j \in \mathbb{N}_0$ not be a multiple of $10$, while $b,c \in 
\mathbb{Z}^+$, and assume radix-$10$ (as usual).
Then, the congruence speed of the tetration $^{b}a \coloneqq j\cdot 10^c$  (i.e., the last digit of $a$ is equal to $0$), is given by the difference between the number of the rightmost (trailing) $0$'s at height $b$ and the number of the rightmost $0$'s at height $b-1$.

Given the fact that $j \not\equiv 0 \pmod{10}$ trivially implies $j^{\left(^{b-1}\left({j\cdot 10^c}\right)\right)} \not\equiv 0 \pmod{10}$, the number of the rightmost $0$'s of \hspace{1mm}$^{b}\left({j \cdot 10^c}\right)$ is the same as the ending $0$'s of $10^{c \cdot \left(^{b-1}\left({j \cdot 10^c}\right)\right)}$ since $^{b}\left({j \cdot 10^c}\right)=\hspace{0.5mm}j^{\left(^{b-1}\left({j \cdot 10^c}\right)\right)} \cdot 10^{\left(c \cdot \left(^{b-1}\left({j \cdot 10^c}\right)\right) \right)}$.

Thus, the tetration $^{b}(j \cdot 10^c )$  is always congruent to $0$ modulo $10^{\left(c \cdot \left(^{b-1}(j \cdot 10^c)\right)\right)}$, whereas $^{b}(j \cdot 10^c) \not\equiv 0 \pmod {10^{\left(1+c \cdot \left(^{b-1}(j \cdot 10^c)\right)\right)}}$, and it follows that ${^b}(j \cdot 10^c)$ has a total of (exactly) $c \cdot \left( ^{b-1}(j \cdot 10^c)\right)$ trailing zeros on the right side.

Hence, as long as $j$, $b-1$, and $c$ are positive integers,
\begin{equation} \label{eq1}
V(j \cdot 10^c,b)=c \cdot \left({^{b-1}(j \cdot 10^c)}\right) - c \cdot \left({^{b-2}(j \cdot 10^c)}\right)  =c \cdot \left(^{b-1}(j \cdot 10^c )- \hspace{0.5mm} ^{b-2}(j \cdot 10^c ) \right)
\end{equation}
(e.g., if $a=300=3 \cdot 10^2$ and $b=3$, then $V(3 \cdot 10^2,3)=2 \cdot (300^{300}-300)=2 \cdot 3^{300}-600$).

Now, we observe that $0^0=1$ implies that 
$^{b}0=1$ if and only if $b$ is even, whereas $^{b}0=0$ for any odd value of $b$.

Consequently, $V(0,b)=0$ for every positive integer $b$.

Therefore, we can conclude that
\begin{equation}\label{eq2}
a \coloneqq j \cdot 10^c \Rightarrow V(a,b)=
\begin{cases}
\hspace{27.4mm} 0 \hspace{27.9mm}\textnormal{iff} \hspace{2mm} j=0\\
\hspace{27.5mm} c \hspace{28mm} \textnormal{iff} \hspace{2mm} j \neq 0 \wedge b=1\\
c \cdot \left( ^{b-1}(j \cdot 10^c)- \hspace{0.5mm} ^{b-2}(j \cdot 10^c) \right) \quad \forall j, b-1, c \in \mathbb{Z}^+
\end{cases},
\end{equation}
holds as long as $j \in \mathbb{N}_0$ is not a multiple of $10$, while $b$ and $c$ are two (strictly) positive integers.

%%%%%%%%%%%%%%%%%%%%%%%%%
%%%%%% CONCLUSION %%%%%%%
%%%%%%%%%%%%%%%%%%%%%%%%%

\section{Conclusion} \label{sec:Concl}

Since we can recursively extend integer tetration so that $^{-1}a$ is well-defined (see Section~\ref{sec:Intr}), Equation (\ref{eq2}), which gives the value of $V(a,b)$ for all the tetration bases that have $0$ as their rightmost digit, can be compactly rewritten (for all $b \in \mathbb{Z}^+$) as 
$$ V(j \cdot 10^c,b)=
\begin{cases}
\hspace{28mm} 0 \hspace{28mm}\textnormal{iff} \hspace{2mm} j=0\\
c \cdot \left( ^{b-1}(j \cdot 10^c)- \hspace{0.5mm} ^{b-2}(j \cdot 10^c) \right) \quad  \textnormal{otherwise}
\end{cases}, $$
and this concludes the present note, together with the description of the congruence speed of the tetration bases that we did not include in our previous papers.

\makeatletter
\renewcommand{\@biblabel}[1]{[#1]\hfill}
\makeatother

\end{document}